\documentclass[review]{elsarticle}
\usepackage{array}
\usepackage{amsthm}
\usepackage{lineno}
\usepackage{textcomp}
\usepackage{amsmath,array,booktabs}
\usepackage{amssymb}
\usepackage{graphics}
\usepackage{amsthm}
\usepackage{CJK}
\usepackage{soul}

\usepackage{hyperref} 
\hypersetup{
	hidelinks,
	hypertex=true,
	colorlinks=true,
	linkcolor=blue,
	citecolor=blue
}
\usepackage{xcolor,cleveref}
\biboptions{numbers,sort&compress}
\usepackage[]{caption2} 
\usepackage[top=1.6cm,bottom=1.6cm,left=1.5cm,right=1.5cm]{geometry}
\makeatletter
\@addtoreset{equation}{section}
\makeatother

\modulolinenumbers[3]
\journal{Journal of \LaTeX\ Templates}

\begin{document}

	\begin{frontmatter}
		
\title{Minimal-norm solution to the Fredholm integral equations of the first kind via the \textit{H}-$H_{K}$ formulation}
		%\author{Renjun Qiu, Wei Qu\corref{mycorrespondingauthor}}
		\author[author1]{Renjun Qiu}
%		\author[author2]{Qian Shang}
		%\author{Renjun Qiu\corref{mycorrespondingauthor}{}}
		\author[author1]{Ming Xu}
		\author[author2]{Wei Qu\corref{mycorrespondingauthor}}
		\cortext[mycorrespondingauthor]{Corresponding author}
		%\fnref{myfootnote}
		\address[author1]{Computer and Information Engineering College, Guizhou University of Commerce,  Guiyang, 550025, China}
		\address[author2]{College of Sciences, China Jiliang University, Zhejiang, 310018, China}
%		\address[author1]{Computer and Information Engineering College, Guizhou University of Commerce,  Guiyang, 550025, China}
%		\address[author2]{Experimental Teaching Center, Guizhou University of Commerce,  Guiyang, 550025, China}
		\fntext[myfootnote]{qiurenjun@whu.edu.cn}
		\fntext[myfootnote]{xminglw@163.com}
		\fntext[myfootnote]{quwei2math@qq.com}
		\begin{abstract}
			%本文研究了第一类Fredholm积分方程的精确最小范数解。具体对于再生核Hilbert空间和平方可积空间，本文已呈现出积分方程的精确最小范数解，它是矩阵方程特解的延拓。
The Fredholm integral equations of the first kind is a typical ill-posed problem, so that it is usually difficult to obtain its analytical minimal-norm solution. This paper gives a closed-form minimal-norm solution for the degenerate kernel equations based on the \textit{H}-$H_{K}$ formulation. Furthermore, it has been shown that the structure of solutions to degenerate kernel equations and matrix equations are consistent. Subsequently, the obtained results are extended to non-degenerate integral equations. Finally, the validity and applicability of the proposed method are demonstrated by some examples. 
		\end{abstract}
		
		\begin{keyword}
			Fredholm integral equations of the first kind\sep Ill-posed problem\sep  Minimal-norm solution \sep Reproducing kernel Hilbert space\sep \textit{H}-$H_{K}$ formulation
		\end{keyword}
	\end{frontmatter}

\section{Introduction}
Consider the Fredholm integral equations of the first kind
\begin{equation}\label{OriginalEq}
	\int_{E} k(x,t)u(t) d t=f(x),~x\in D,
\end{equation}
where the function $u(t)$ is unknown and to be determined, $k(x,t)$ and $f(x)$ are given. Usually in the sense of $L^{2}-$norm, Eq.~\eqref{OriginalEq} is ill-posed, in other words, it does not satisfy one of the Hadamard conditions at least
\begin{itemize}
	\item Eq.~\eqref{OriginalEq} has a solution $u(t)$,
	\item Solution $u(t)$ is unique, and
	\item Continuous dependence of $u(t)$ on $f(x)$.
\end{itemize}
Therefore it is difficult to find feasible numerical solutions, let alone feasible analytical solutions. To address this problem, from another perspective, literatures \cite{golbabai2008modified,Wazwaz2011regularization,golbabai2009solution,alturk2019use,alturk2016regularization,yuldashev2023nonlinear} have given the analytical solutions for a special type of the integral kernel, i.e., degenerate kernel 
\begin{equation}\label{DegenerateKernel}
	k(x,t)=\sum\limits_{i=1}^{n} g_{i}(x) h_{i}(t), n\geq 1.
\end{equation}

In reality, as long as one obtains analytical solutions of the degenerate kernel Fredholm integral equation (DKFIE), then solutions of the non-degenerate kernel Fredholm integral equation (NDKFIE) can be  approximated well \cite{qiu2022solving,qiu2022best}. On account of this reason, exploring analytical solutions of the DKFIE~\eqref{OperatorEq} is very necessary to a NDKFIE~\eqref{OriginalEq}, in particular, the analytical minimal-norm solution.

Without loss of generality, we assume that $\{g_{i}(x)\}_{i=1}^{n}$ and $\{h_{i}(t)\}_{i=1}^{n}$ are linearly independent functions in this paper. In despite of literatures have appled the modified homotopy (perturbation) methods \cite{golbabai2008modified,alturk2016regularization,golbabai2009solution} or regularization methods \cite{Wazwaz2011regularization,alturk2019use,yuldashev2023nonlinear} to obtain many solutions of the DKFIE, yet not the minimal-norm solution. The motivation of this paper is to excavate a closed-form representation of the minimal-norm solution, furthermore, to explore the structure of solutions to a DKFIE, imitating a linear matrix equation. 

To achieve this target, based on the \textit{H}-$H_{k}$ formulation given in \cite{qian2020reproducing}, the minimal-norm solution of the DKFIE has been obtained by us, which is a crucial part of the structure of solutions. Next, we give the structure of solutions for a DKFIE. Similar to a linear matrix equation, any solution in which can be decomposed into two parts, i.e., one is in the null space and the other is in the null complement space. Finally, the obtained results are extended form DKFIE to NDKFIE.

\section{Preliminaries}
For the convenience of writing, Eq.~\eqref{OriginalEq} can be abbreviated as
\begin{equation}\label{OperatorEq}
	L(u)(x):=\langle k_{x}, u\rangle_{L^{2}(E)}=f(x),
\end{equation} 
where $k_{x}(t):=k(x,t)=G^{T}(x)H(t)$,~$G^{T}(x)=(g_{1}(x),\cdots,g_{n}(x))$,~$H^{T}(t)=(h_{1}(t),\cdots,h_{n}(t))$, i.e., Eq.~\eqref{OperatorEq} is a DKFIE.
\subsection{Minimal-norm solution}
Let $ N(L)$ and $R(L)$ be denoted the null space and range space of $L$, $ N(L)^{\perp} $ be the null complementary space of $ N(L)$, and  $ P_{N(L)^{\perp}} $ be  the orthogonal projection from  $L^{2}(E)$ onto $N(L)^{\perp}$. Let $f(x)\in L^{2}(D)$, a function $u \in L^{2}(E)$ of DKFIE~\eqref{OperatorEq} is called a \textit{least-squares solution} \cite{pes2023projection}, if 
\begin{equation}
	\|Lu-f\|_{L^{2}(D)}=\inf \left\{\|Lv-f\|_{L^{2}(D)}: v \in L^{2}(E)\right\}.
\end{equation}
An element $u^{\dagger}\in S$ is called the \textit{minimal-norm solution} of the DKFIE~\eqref{OperatorEq}, if $\|u^{\dagger}\|_{L^{2}(E)}=\inf\limits _{u \in S}\|u\|_{L^{2}(E)},$
here $S$ is the set composed of least-squares solutions. Meanwhile, $L^{\dagger}$ is referred to as the \textit{Moore-Penrose inverse operator}, if $L^{\dagger} f:=u^{\dagger}$. 

For the minimal-norm solution of the NDKFIE~\eqref{OriginalEq}, some interesting discussions and remarks may be found in \cite{bechouat2023collocation,bechouat2023numerical}.
%\newtheorem{Definition}{Definition}
%\begin{Definition}
%	\textit{Eq.~\eqref{OriginalEq} is called} a solvable equation, when it satisfies $g\in R(L)$; Eq.~\eqref{OriginalEq} is called a solvable projection equation, when it satisfies $P_{\overline{R(L)}}g\in R(L).$  
%\end{Definition}

\newtheorem{Proposition}{Proposition}
%\begin{Proposition}\label{Prop1}Let $f(x)\in L^{2}(D)$, then
%	$u^{\dagger}$ existss if and only if $ P_{R(L)^{\perp}}u $.
%\end{Proposition}	
\subsection{\textit{H}-$H_{K}$ formulation}
In this subsection, the \textit{H}-$H_{K}$ formulation originated by \cite{qian2020reproducing} is introduced to study a DKFIE~\eqref{OperatorEq}, which can establish an isometric isomorphism between $N(L)^{\perp}$ and $R(L)$.

Firstly, $N(L)^{\perp}$ can be described accurately, which determines the representations of all solutions of the DKFIE~\eqref{OperatorEq}.
\newtheorem{Lemma}{Lemma}%[section 2]
\begin{Lemma}\label{lemma2} $N(L)^{\perp}$ can be represented as
	\begin{equation}\label{ZeroSpace}
		N(L)^{\perp}=span\left\{h_{1}(t),\cdots,h_{n}(t),~t \in E\right\}.
	\end{equation}
	%\noindent where $\eta_{x}(t)=\int_{E} h_{x}(s)Q(t,s) ds$, and $Q(t,s)$ is the reproducing kernel of RKHS $H_{Q}$.
\end{Lemma}
\noindent $\textbf{Proof.}$ Let $u\in N(L)$, since $\{g_{i}(x)\}_{i=1}^{n}$ are linearly independent functions in DKFIE~\eqref{OperatorEq}, then 
\begin{equation}\label{NullSpace}
	\int_{E}h_{i}(t)u(t) d t=0,~i=1,\cdots,n.
\end{equation}
That is, $h_{i}(t)\in N(L)^{\perp}$, thereby $span\left\{h_{1}(t),\cdots,h_{n}(t),~t \in E\right\}\subseteq N(L)^{\perp}.$ 

Let $u\in span\left\{h_{1}(t),\cdots,h_{n}(t),~t \in E\right\}^{\perp}$, then Eq.~\eqref{NullSpace} holds. Consequently, we have
\begin{equation*}
	L(u(t))(x)=\sum\limits_{i=1}^{n} g_{i}(x)\int_{E}h_{i}(t)u(t) d t=0,
\end{equation*}
that is, $u\in N(L)$, $span\left\{h_{1}(t),\cdots,h_{n}(t),~t \in E\right\}^{\perp}\subseteq N(L)$. Hence~$N(L)^{\perp}\subseteq span\left\{h_{1}(t),\cdots,h_{n}(t),~t \in E\right\}$.
~~\hfill $\square$

Secondly, the range space $R(L)$ designated a specific norm becomes a RKHS. 

For any given $f_{1}, f_{2}\in R(L)$, an inner product can be seen in  \cite{qiu2022solving,qiu2022best,qian2020reproducing} defined by
\begin{equation}\label{InnerProduct}
	\left\langle f_{1}, f_{2}\right\rangle_{R(L)}\!:=\left\langle L^{\dagger}f_{1}, L^{\dagger}f_{2}\right\rangle_{L^{2}(E)}.
\end{equation}
\begin{Lemma}\label{Lemma1}
 Under the inner product \eqref{InnerProduct},	$R(L)$	
is a RKHS with reproducing kernel defined by
	\begin{equation}\label{ReproducingKernel}
		K(x,x^{\prime})=G^{T}(x)HG(x^{\prime}),~x,x^{\prime} \in D.
	\end{equation}
where $\displaystyle{H:=H_{n\times n}=[h_{ij}]=[\int_{E}h_{i}(t)h_{j}(t) d t]}$. Moreover, there exists $x_{1},\cdots,x_{n}\in D$ such that
\begin{equation}\label{RangeSpace}
	R(L)=span\left\{K(x,x_{1}),\cdots,K(x,x_{n}),~x \in D\right\}.
\end{equation}	 
\end{Lemma}
\noindent $\textbf{Proof.}$ For given $x^{\prime}\in D$, we have
\begin{equation*}
	 L(G^{T}(x^{\prime})H(t))(x)=	\int_{E}G^{T}(x^{\prime})H(t)H^{T}(t)G(x)dt=G^{T}(x)HG(x^{\prime})=K(x,x^{\prime}).
\end{equation*}
Note that $K_{x^{\prime}}(x):=K(x,x^{\prime})$, then $L(G^{T}(x^{\prime})H(t))=K_{x^{\prime}}$. By lemma \ref{lemma2}, we get 
\begin{equation*}
	G^{T}(x^{\prime})H(t)\in N(L)^{^{\perp}}.
\end{equation*}

For any $f(x)\in R(L)$, there exists $C\in R^{n}$ such that $L(C^{T}H(t))(x)=f(x)$, i.e., 
\begin{equation*}
	G^{T}(x)HC=f(x),
\end{equation*}
then we have
\begin{equation*}
	\left\langle f, K_{x}\right\rangle_{R(L)}=\left\langle C^{T}H(t), G^{T}(x)H(t)\right\rangle_{L^{2}(E)}=G^{T}(x)HC=f(x).
\end{equation*}
That is to say, $K$ defined by Eq.~\eqref{ReproducingKernel} is a reproducing kernel in $R(L)$, and $R(L)$	
is a RKHS. 

Since
%\begin{equation*}
%	R(L)=span\left\{g_{1}(x),\cdots,g_{n}(x),~x \in D\right\}
%\end{equation*}
$\{g_{i}(x)\}_{i=1}^{n}$ are linearly independent, then there exists $x_{1},\cdots,x_{n}\in D$ such that
\begin{equation*}
	(g_{1}(x_{1}),\cdots,g_{n}(x_{n}))
\end{equation*}
becomes an invertible matrix, then $R(L)=span\left\{K(x,x_{1}),\cdots,K(x,x_{n}),~x \in D\right\}$ by Eq.~\eqref{ReproducingKernel}.~~\hfill $\square$

%Based on Eq.~\eqref{ReproducingKernel}, we obtain 
%\begin{equation}\label{IsometricIsomorphism}
%	L\eta_{x}=K_{x},~Lh_{x}=K_{x}.
%\end{equation}
%Note that $\eta_{x}, h_{x}\in N(L)^{\perp}$, $K_{x}\in H_{K}(=R(L))$, then $L^{\dagger}K_{x}=\eta_{x}$ or $L^{\dagger}K_{x}=h_{x}$, thereby $L^{\dagger}$ is a one-to-one bijective mapping from $g$ to $L^{\dagger}g$, i.e., $N(L)^{\perp}$ is isometric isomorphism of $H_{K}$. This isomorphic relationship \eqref{IsometricIsomorphism} is called \textit{the H-$H_{K}$ formulation} in \cite{qian2020reproducing}. More detailed properties and applications of the \textit{H}-$H_{K}$ formulation could be found in \cite{qian2020reproducing}.

\section{Structure of solutions}
Note that $\displaystyle{F:=F_{n\times1}=[f_{k}]=[\int_{E} f(x)h_{k}(t) d t]}$,~$\displaystyle{A:=A_{n\times n}=[a_{ij}]=[\int_{E}h_{i}(t)g_{j}(t) d t]}$. Based on these notations, the following analytic solution has been obtained for literature \cite{alturk2019use}.

\begin{Proposition}\label{Prop2}
Let $A$ be invertible and $f(x)=(A^{-1}F)^{T} G(x)$, then DKFIE~\eqref{OperatorEq} has a solution 
\begin{equation}\label{BeforeSolution}
	u(x)=((A^{-1})^2F)^{T}G(x),
\end{equation}
\begin{equation}\label{BeforeSolutionNorm}
	\| u\|^{2}_{L^{2}}=((A^{-1})^2F)^{T}G((A^{-1})^2F),
\end{equation}
where $\displaystyle{G:=G_{n\times n}=[g_{ij}]=[\int_{E}g_{i}(x)g_{j}(x) dx]}$.
\end{Proposition}
This analytical solution \eqref{BeforeSolution} is available, which is not our target. Our targets are twofold, one is to give the minimal-norm solution and structure of solutions of the DKFIE, the other is to extend the results obtained to a NDKFIE. To achieve these targets, we assume that functions $\left\lbrace \varphi_{i}(t)\right\rbrace_{i=1}^{\infty}$ are a basis in $N(L)$.
%
%Since the latter is necessary to explore the structure of solutions, and the \textit{H}-$H_{K}$ formulation which are summarized in the following theorem.
\newtheorem{Theorem}{Theorem}
\begin{Theorem}\label{Theorem1} 
Let $A$ be invertible and $f(x)=(A^{-1}F)^{T} G(x)$, then DKFIE~\eqref{OperatorEq} has the minimal-norm solution 
\begin{equation}\label{Minimal-normSolution}
	u^{\dagger}(t)=(H^{-1}A^{-1}F)^{T}H(t),
\end{equation}	
\begin{equation}\label{Minimal-normSolutionNorm}
	\| u^{\dagger}\|^{2}_{L^{2}}=(A^{-1}F)^{T}H^{-1}(A^{-1}F).
\end{equation}
Moreover, the structure of solutions can be represented as
\begin{equation}\label{StructureSolutions}
	u(t)=u^{\dagger}(t)+\sum\limits_{i=1}^{\infty} c_{i}\varphi_{i}(t),
\end{equation}
where the real sequences $\left\lbrace c_{i}\right\rbrace_{i=1}^{\infty}\in l^{2} $, i.e., $\sum\limits_{i=1}^{\infty}c^{2}_{i}<+\infty$.
\end{Theorem}
\noindent$\mathbf{Proof.}$ By $f(x)=(A^{-1}F)^{T} G(x)$, DKFIE~\eqref{OperatorEq} is solvable, i.e., $f(x)\in R(L)$. By lemma \ref{Lemma1}, there exists $x_{1},\cdots,x_{n}\in D$, as well as matrix $U$ such that
\begin{equation*}
	G(x)=UK_{X}(x),
\end{equation*}
where $K^{T}_{X}(x)=(K(x,x_{1}),\cdots,K(x,x_{n}))$. Based on Eq.~\eqref{ReproducingKernel}, we get 
\begin{equation*}
	K_{X}(x)=G_{X}^{T}HG(x),
\end{equation*}
where $G_{X}:=[g_{i}(x_{j})]$. Hence, we get $U=H^{-1}(G_{X}^{T})^{-1}$ and 
\begin{equation*}\label{AnotherStyle}
	f(x)=F^{T}(A^{-1})^{T}H^{-1}(G_{X}^{T})^{-1}K_{X}(x).
\end{equation*}
Based on the \textit{H}-$H_{K}$ formulation, for given $ x_{i}\in D,~1\leq i\leq n$, we have
\begin{equation*}
	L[G^{T}(x_{i})H(t)]=K(x,x_{i}),~
%\end{equation*}
%\begin{equation*}
	L^{\dagger}K_{X}(x)=G_{X}^{T}H(t).
\end{equation*}
According to the \textit{H}-$H_{K}$ formulation again, we get 
\begin{equation*}
	u^{\dagger}(t)=L^{\dagger}f(x)=F^{T}(A^{-1})^{T}H^{-1}(G_{X}^{T})^{-1}L^{\dagger}K_{X}(x)=(H^{-1}A^{-1}F)^{T}H(t).
\end{equation*}
\begin{equation*}
	\|u^{\dagger}\|^{2}_{L^{2}}=\int_{E}(H^{-1}A^{-1}F)^{T}H(t)H^{T}(t)H^{-1}A^{-1}Fdt=(A^{-1}F)^{T}H^{-1}(A^{-1}F).
\end{equation*}

Since $\left\lbrace \varphi_{i}(t)\right\rbrace_{i=1}^{\infty}$ is a basis in $N(L)$, as well as $\left\lbrace h_{i}(t)\right\rbrace_{i=1}^{n}$ is a basis in $N(L)^{\perp}$, then
\begin{equation}
	\left\lbrace h_{i}(t)\right\rbrace_{i=1}^{n}\cup\left\lbrace \varphi_{i}(t)\right\rbrace_{i=1}^{\infty}
\end{equation}
is a complete basis in $L^{2}(E)$, then any solution $u(t)$ can be decomposed as
\begin{equation*}\label{StructureSolutions}
	u(t)=u^{\dagger}(t)+\sum\limits_{i=1}^{\infty} c_{i}\varphi_{i}(t),
\end{equation*}
where  the real sequences $\left\lbrace c_{i}\right\rbrace_{i=1}^{\infty}\in l^{2} $.  \hfill $\square$

Actually, according to the process of proof, $A^{-1}F$ can be substituted fully by a given vector $C\in R^{n} $, which can ensure the existence of a solution to DKFIE. Namely, we no longer need $A$ invertible in the present paper, see example 3.
\newtheorem{Corollary}{Corollary}
\begin{Corollary}\label{Corollary1}
Let $f(x)=C^{T}G(x)$ for a given vector $ C\in R^{n} $ in DKFIE~\eqref{OperatorEq}, then
\begin{equation}\label{GeneralSlution}
	u^{\dagger}(t)=C^{T}H^{-1}H(t),
\end{equation}
\begin{equation}\label{GeneralSlutionNorm}
	\| u^{\dagger}\|^{2}_{L^{2}}=C^{T}H^{-1}C.
\end{equation}
\end{Corollary}
\begin{Corollary}\label{Corollary2}
Let $C=A^{-1}F$ in corollary \ref{Corollary1} and $G(x)=KH(x)$ for a given invertible matrix $ K$, then \eqref{BeforeSolution} and \eqref{Minimal-normSolution}, as well as \eqref{BeforeSolutionNorm} and \eqref{Minimal-normSolutionNorm} are consistent.
\end{Corollary}

Finally, for a general NDKFIE, the minimal-norm solution can be obtained by discussing analogously under the \textit{H}-$H_{K}$ formulation. In fact, assuming that
\begin{equation}\label{Non-DegenerateKernel}
k(x,t)=\sum\limits_{i=1}^{\infty} g_{i}(x) h_{i}(t),
%\lim_{n\to \infty}\int_{D} \int_{E}\left(k(x,t)-\sum\limits_{k=1}^{n} g_{k}(x) h_{k}(t)\right)^{2} dtdx=0,
\end{equation}
\begin{equation}\label{Non-DegenerateTerm}
%\lim_{n\to \infty}\int_{D} \left(f(x)-f_{n}(x)\right)^{2}dx=0,
f(x)=\sum\limits_{i=1}^{\infty} c_{i}g_{i}(x).
\end{equation}
%where the real sequences $\left\lbrace c_{i}\right\rbrace_{i=1}^{\infty}\in l^{2} $ again. 
\begin{Theorem}\label{Corollary3} 
Under the assumptions of Eq.~\eqref{Non-DegenerateKernel} and Eq.~\eqref{Non-DegenerateTerm}, then NDKFIE ~\eqref{OriginalEq} has the minimal-norm solution 
\begin{equation}\label{Non-DegenerateSolution}
%	\lim_{n\to \infty}\int_{E} \left(u^{\dagger}(t)-u^{\dagger}_{n}(t)\right)^{2}dt=0,
u^{\dagger}(t)=C^{T}H^{-1}H(t),~	
\end{equation}
\begin{equation}\label{Non-DegenerateNorm}
	\| u^{\dagger}\|^{2}_{L^{2}}=C^{T}H^{-1}C,
\end{equation}
where $C=[c_{i}]$ and $H(t)$ are infinite dimensional column vectors, and $H$ is an infinite dimensional matrix. 
\end{Theorem}
This theorem extends corollary \ref{Corollary1} from DKFIE to NDKFIE, which is the desired result, see example 5. The proof of this theorem can be  modified directly on the basis of theorem \ref{Theorem1}, so it will not be elaborated here.
%	\begin{equation}\label{GeneralSlution}
%		u^{\dagger}(t)=(H^{-1}C)\cdotp H(t),
%	\end{equation}
%	\begin{equation}\label{GeneralNorm}
%		\| u^{\dagger}\|^{2}_{L^{2}}=C^{T}H^{-1}C.
%	\end{equation}

\newtheorem{Remark}{Remark}
\begin{Remark}
For DKFIEs, literatures \cite{golbabai2008modified,alturk2016regularization,golbabai2009solution,Wazwaz2011regularization,alturk2019use,yuldashev2023nonlinear} can indeed provide analytical solutions, which have the following form
\begin{equation*}
	u(t)=\sum\limits_{i=1}^{n} c_{i}g_{i}(t),
\end{equation*}
where $c_{i},1\leq i\leq n$ are determined by the modified homotopy methods  or regularization methods. Note that existing methods require the following conditions  
\begin{equation}\label{OrthogonalConditions}
	\int_{E}h_{i}(t)g_{j}(t) d t\neq0,~1\leq i,j\leq n,
\end{equation}
to be met, that is, $g_{j}(t)$ cannot be orthogonal to $h_{i}(t)$. Specifically, A is invertible in proposition \ref{Prop2}. It is essentially different from minimal-norm solutions \eqref{GeneralSlution} \eqref{Non-DegenerateSolution} presented in this paper, as they are represented by functions $h_{i}(t),1\leq i\leq n$, which no longer need to meet  \eqref{OrthogonalConditions}.
\end{Remark}	

%In summary, there are two typical analytical solutions, when the right-hand term $g(x)$ has been obtained analytically. One is that the solution can be assumed $f(t)=\sum\limits_{i=1}^{m} c_{i}a_{i}(t)$, where $\{a_{i}\}_{i=1}^{m}$ comes from Eq.~\eqref{DegenerateKernel}.   in, if $g(x)=A_{m}^{T}(x)(A^{-1})^{T}F$, then 
%\begin{equation}\label{FirstSolution}
%	f(t)=A_{m}^{T}(x)(A^{-2})^{T}F,
%\end{equation}
%where $A=\left[\int_{E} b_{i}(t) a_{j}(t) dt\right]$, $F=\left(\int_{E} g(t) b_{1}(t) dt,\cdots,\int_{E} g(t) b_{m}(t) dt\right)^{T}$, and $A_{m}(x)=\left( a_{1}(x),\cdots,a_{m}(x)\right)^{T}$. Another is that the solution can be assumed  $f(t)=\sum\limits_{i=1}^{m} c_{i}b_{i}(t)$, where $\{b_{i}\}_{i=1}^{m}$ also comes from Eq.~\eqref{DegenerateKernel}. In \cite{alturk2019use}, if $g(x)=A_{m}^{T}(x)C$ for a certain column vector $C$, then 
%\begin{equation}\label{SecondSolution}
%	f(t)=B_{m}^{T}(t)B^{-1}C,
%\end{equation}
%where $B=\left[\int_{E} b_{i}(t) b_{j}(t) dt\right]_{m\times m} $, and $B_{m}(t)=\left(b_{1}(t),\cdots, b_{m}(t)\right)^{T}$.

\section{Illustrative Examples}
To start with, we compare solutions of DKFIEs obtained by our methods with solutions by \cite{golbabai2008modified,golbabai2009solution,alturk2019use,bechouat2023numerical} for the same examples. Next, we provide an example of the DKFIE in which $A$ is non-invertible, its minimal-norm solution can also be obtained by the proposed method. Finally we discussed the minimal-norm solution of a NDKFIE derived from an one-dimensional backward heat conduction problem. 

\noindent \textbf{Example 1.} Calculate the minimal-norm solution of the DKFIE 
\begin{equation}\label{ExampleEquation3}
	\int_{0}^{1} (te^{x}+1)u(t) d t=\frac{1}{3}e^{x}+\frac{1}{2}.
\end{equation}	
By direct calculation 
\begin{equation*}
	%	A=\left(\begin{array}{cc}
		%		5/3 & 5/4 \\
		%		5/4 & 1
		%	\end{array}\right),
	H=\left[ \begin{array}{cc}
		1 & 1/2 \\
		1/2 & 1/3
	\end{array}\right],~
	H^{-1}=\left[ \begin{array}{cc}
		4 & -6 \\
		-6 & 12
	\end{array}\right] ,
	~C=\left[ \begin{array}{c}
		1/2 \\
		1/3
	\end{array}\right] ,~
	%	H=\left(\begin{array}{cc}
		%		1/3 & 1/4 \\
		%		1/4 & 1/5
		%	\end{array}\right),~	
	H(t)=\left[ \begin{array}{c}
		1 \\
		t
	\end{array}\right] 
\end{equation*}	
according to corollary \ref{Corollary1}, then
\begin{equation*}
	u^{\dagger}(t)=C^{T}H^{-1}H(t)=t,~\| u^{\dagger}\|_{L^{2}}=\sqrt{3}/3.
\end{equation*}
There has another solution 
\begin{equation*}
	u(t)=\frac{1}{18-6e}e^{t}+\frac{5-2e}{9-3e}
\end{equation*}
which is obtained from \cite{golbabai2008modified}. 

Similarly, for example 2 in \cite{golbabai2009solution}, by corollary \ref{Corollary1}, we have 
\begin{equation*}
	e^{u^{\dagger}(t)}=(18-6e)t+4e-10~(u^{\dagger}(t)=ln[(18-6e)t+4e-10]).
\end{equation*}
There has also another solution 
\begin{equation*}
		e^{u(t)}=e^{t}.
\end{equation*}

\noindent \textbf{Example 2.} Calculate the minimal-norm solution of the DKFIE
\begin{equation}\label{ExampleEquation1}
		\int_{0}^{\pi} \cos x\sin t~u(t) d t=\frac{\pi}{2} \cos x.
\end{equation}
Obviously, $f(x)=\frac{\pi}{2}\cos x$,~$g(x)=\cos x$,~ $h(t)=\sin t$. By corollary \ref{Corollary1}, then
\begin{equation*}
	u^{\dagger}(t)=\sin t,~\| u^{\dagger}\|_{L^{2}}=\sqrt{\pi/2}.
\end{equation*}	

Since $f(t)\perp h(t)$ in $L^{2}([0,\pi])$, in this case,  an auxiliary function $\phi_{\beta}(t)$ needs to be introduced in \cite{alturk2019use} such that
\begin{equation*}
	u(t)=\frac{\pi\cos t\phi_{\beta}(t)}{2\int_{0}^{\pi}\cos t\sin t\phi_{\beta}(t)dt}
\end{equation*}
becomes a class of solutions of Eq.~\eqref{ExampleEquation1}. Because
\begin{equation}\label{CauchyEq1}
	(\int_{0}^{\pi}\cos t\sin t\phi_{\beta}(t)dt)^{2}\leq\int_{0}^{\pi}\cos^{2}(t)\phi_{\beta}^{2}(t)dt\int_{0}^{\pi}\sin^{2}(t)dt=\frac{\pi}{2}\int_{0}^{\pi}\cos^{2}(t)\phi_{\beta}^{2}(t)dt
\end{equation}
\begin{equation}\label{CauchyEq2}
\int_{0}^{\pi}u^{2}(t)dt\geq \frac{\pi^{2}}{4}\frac{\int_{0}^{\pi}\cos^{2}(t)\phi_{\beta}^{2}(t)dt}{\frac{\pi}{2}\int_{0}^{\pi}\cos^{2}(t)\phi_{\beta}^{2}(t)dt}=\frac{\pi}{2},
\end{equation}
then $\| u^{\dagger}\|_{L^{2}}\leq\| u\|_{L^{2}}$ for any given $\phi_{\beta}(t)$. 

Similarly, example 2 in \cite{alturk2019use} can also be calculated and compared consequentially.

%\noindent \textbf{Example 2.} Calculate the minimal-norm solution of the Fredholm integral equation
%\begin{equation}\label{ExampleEquation2}
%	\int_{0}^{1} e^{x-2t}u(t) d t=e^{x}.
%\end{equation}	
%Obviously, $f(x)=g(x)=e^{x}$,~$h(t)=e^{-2t}$. By \eqref{GeneralSlution} and \eqref{GeneralNorm} in corollary \ref{Corollary1}, we obtain the minimal-norm solution 
%\begin{equation*}
%	u^{\dagger}(t)=\frac{4}{1-e^{-4}}e^{-2t},~\| u^{\dagger}\|_{L^{2}}=\sqrt{\frac{4}{1-e^{-4}}}.
%\end{equation*}	
%In \cite{alturk2019use}, an auxiliary function $\phi_{\beta}(x)$ is introduced such that
%\begin{equation*}
%	u(t)=\frac{e^{t}\phi_{\beta}(t)}{\int_{0}^{1}e^{t-2t}\phi_{\beta}(t)dt}
%\end{equation*}
%becomes a solution of Eq.~\eqref{ExampleEquation1}. We imitate \eqref{CauchyEq1} and \eqref{CauchyEq2}, we also have
%\begin{equation*}
%	u^{\dagger}\|_{L^{2}}\leq \| u\|_{L^{2}}
%\end{equation*}
%for any given auxiliary function $\phi_{\beta}(x)$.

\noindent \textbf{Example 3.} Calculate the minimal-norm solution of the DKFIE
\begin{equation}\label{ExampleEquation3}
	\int_{0}^{1} 5(xt+x^{2}t^{2})u(t) d t=x+6x^{2}.
\end{equation}	
By direct calculation 
\begin{equation*}
%	A=\left(\begin{array}{cc}
%		5/3 & 5/4 \\
%		5/4 & 1
%	\end{array}\right),
A^{-1}=\left[ \begin{array}{cc}
		48/5 & -12 \\
		-12 & 16
	\end{array}\right],~
H^{-1}=\left[ \begin{array}{cc}
	48 & -60 \\
	-60 & 80
\end{array}\right] ,
~F=\left[ \begin{array}{c}
		11/6 \\
		29/20
	\end{array}\right] ,~
%	H=\left(\begin{array}{cc}
%		1/3 & 1/4 \\
%		1/4 & 1/5
%	\end{array}\right),~	
H(t)=\left[ \begin{array}{c}
	t \\
	t^{2}
\end{array}\right] 
\end{equation*}	
according to \eqref{Minimal-normSolution} in theorem \ref{Theorem1}, then 
\begin{equation}
	u^{\dagger}(t)=(H^{-1}A^{-1}F)^{T}H(t)=84t^{2}-\frac{312}{5}t.
\end{equation}
	
Since $G(x)=5H(x)$, by corollary \ref{Corollary2}, we have
$u^{\dagger}(t)= u(t),$ where $u(t)$ is a solution obtained by \cite{alturk2019use}. Similarly, for example 4 in \cite{alturk2019use}, since $G(x)=KH(x)$, namely,
\begin{equation*}
	%	A=\left(\begin{array}{cc}
		%		5/3 & 5/4 \\
		%		5/4 & 1
		%	\end{array}\right),
\left[ \begin{array}{c}
		\sin x \\
		\cos x
	\end{array}\right] 
=\left[ \begin{array}{cc}
		0 & -1 \\
		1 & 0
	\end{array}\right] 
\left[ \begin{array}{c}
	\cos x \\
	-\sin x
	\end{array}\right] 
\end{equation*}	
then 
\begin{equation*}
	u^{\dagger}(t)=u(t)=\frac{24}{\pi^{2}-4}(\frac{\pi}{2}\sin t-\cos t)
\end{equation*}
by corollary \ref{Corollary2}, where $u(t)$ is a solution obtained by \cite{alturk2019use}. 

%\noindent \textbf{Example 3.} Calculate the minimal-norm solution of the Fredholm integral equation
%\begin{equation}\label{ExampleEquation4}
%	\int_{0}^{1} (4te^{x}+3)u(t) d t=e^{x+1}.
%\end{equation}	
%By direct calculation 
%\begin{equation*}
%%A=\left(\begin{array}{cc}
%%	4 & 3/2 \\
%%	4(e-1) & 3
%%\end{array}\right),~
%A^{-1}=\displaystyle{\frac{1}{18-6e}}\left(\begin{array}{cc}
%	3 & -3/2 \\
%	4(1-e) & 4
%\end{array}\right), ~
%%	H=\left(\begin{array}{cc}
%%		1/3 & 1/2 \\
%%		1/2 & 1
%%	\end{array}\right),~	
%	H^{-1}=\left(\begin{array}{cc}
%		12 & -6 \\
%		-6 & 4
%	\end{array}\right),~
%F=\left(\begin{array}{c}
%	e \\
%	e(e-1)
%\end{array}\right),~
%	H(t)=\left(\begin{array}{c}
%		t \\
%		1
%	\end{array}\right)
%\end{equation*}	
%according to \eqref{Minimal-normSolution} and \eqref{Minimal-normSolutionNorm}  in theorem \ref{Theorem1}, then 
%\begin{equation*}
%	u^{\dagger}(t)=3e(t-\frac{1}{2}),~\| u^{\dagger}\|_{L^{2}}=\frac{\sqrt{3}e}{2}.
%\end{equation*}
%Since $u(t)=\displaystyle{\frac{e^{2}-e-e^{t+1}}{2(e-3)}}$ in \cite{alturk2019use}, we get 
%\begin{equation*}
%	\| u\|_{L^{2}}=\sqrt{\frac{e^{3}-e^{2}}{24-8e}},~\| u\|_{L^{2}}>\| u^{\dagger}\|_{L^{2}}.
%\end{equation*}

\noindent \textbf{Example 4.} Calculate the minimal-norm solution of the DKFIE
\begin{equation}\label{ExampleEquation4}
	\int_{-\pi/2}^{\pi/2} (\sin x\cos t+1)u(t) d t=\sin x.
\end{equation}	
By direct calculation, $A$ is a non-invertible matrix, i.e., beyond the discussion of \cite{alturk2019use}, and
\begin{equation*}
	H^{-1}=\displaystyle{\frac{2\pi}{\pi^{2}-8}}\left[ \begin{array}{cc}
		\pi/2 & -2 \\
		-2 & \pi
	\end{array}\right] , ~
	C=\left[ \begin{array}{c}
		0 \\
		1
	\end{array}\right] ,~
	H(t)=\left[ \begin{array}{c}
		1 \\
		\cos t
	\end{array}\right]. 
\end{equation*}	
In terms of corollary \ref{Corollary1}, then 
\begin{equation*}
	u^{\dagger}(t)=\frac{2}{\pi^{2}-8}(\pi\cos t-2),~\| u^{\dagger}\|_{L^{2}}=\sqrt{\frac{2\pi}{\pi^{2}-8}}.
\end{equation*}

\noindent \textbf{Example 5.} Calculate the minimal-norm solution of the DKFIE \cite{bechouat2023numerical}
\begin{equation}\label{ExampleEquation4}
	\int_{0}^{1}\int_{0}^{1}e^{\tau^{2}+\eta^{2}+s+t-2} u(s,t) d tds=\frac{1}{4}(e^{-2}-1)^{2}e^{\tau^{2}+\eta^{2}}.
\end{equation}	
Let $f(\tau,\eta)=\frac{1}{4}(e^{-2}-1)^{2}e^{\tau^{2}+\eta^{2}}$, $g(\tau,\eta)=e^{\tau^{2}+\eta^{2}-2}$, then
\begin{equation*}
	f(\tau,\eta)=\frac{e^{2}}{4}(e^{-2}-1)^{2}g(\tau,\eta).
\end{equation*}
By corollary \ref{Corollary1}, we have
\begin{equation*}
	C=\frac{e^{2}}{4}(e^{-2}-1)^{2},~H=\frac{(e^{2}-1)^{2}}{4},~H(s,t)=e^{s+t},
\end{equation*}
then 
\begin{equation*}
u^{\dagger}(s,t)=C^{T}H^{-1}H(s,t)=e^{s+t-2},~\| u^{\dagger}\|_{L^{2}}=\frac{1-e^{-2}}{2}.
\end{equation*}

By the same method proposed by this paper, the minimal-norm solution for example 2 in \cite{bechouat2023numerical} can also be obtained as
\begin{equation*}
	u^{\dagger}(s,t)=(2s-1)(2t-1)^{3}.
\end{equation*}

\noindent \textbf{Example 6.}
Calculate the following initial value problem \cite{du2008approximate}
\begin{equation}\label{one-BHCP}
	\begin{cases}
		u_{t}(x, s) =u_{x x}(x, s), & 0<x<\pi,  \\
		u(0, s)     =u(\pi, s)=0, &  0<s,                 \\
		u(x, 0)     =u_{0}(x), & 0 \leq x \leq  \pi,
	\end{cases}
\end{equation}
where $u_{0}(x)$ needs to be determined.

In terms of separation variable method, there has a formal solution 
\begin{equation}\label{FormalSolution}
	u(x, s)=\sum\limits_{i=1}^{\infty} a_{i} e^{-i^{2}s} \sin ix,~ a_{i}=\frac{2}{\pi} \int_{0}^{\pi} u_{0}(t) \sin it~dt.
\end{equation}
We substitute $a_{i}$ into $u(x, s)$ as $s=1$, then we obtain a NDKFIE as below
\begin{equation}\label{One-BHCP-fredhlom}
	\int_{0}^{\pi} k(x,t) u_{0}(t) \mathrm{d} t=u(x, 1), \quad 0 \leq x \leq \pi.
\end{equation}
Herein $k(x, t)=\frac{2}{\pi}\sum\limits_{i=1}^{\infty} e^{-i^{2}} \sin ix \sin it$ and $u(x,1)=\frac{\sin x}{e}$, which are consistent with \cite{du2008approximate}. 

Let $h_{i}(t)=\sin it,~g_{i}(x)=\frac{2}{\pi}\sin ix, i\in N^{\ast}$, by theorem \ref{Corollary3}, then we have 
\begin{equation*}
	u(x, 1)=\frac{\pi}{2}g_{1}(x),~H^{-1}_{ii}=\frac{2}{\pi},~H^{-1}_{ij}=0,~i\neq j,~i,j\in N^{\ast}.
\end{equation*}
According to \eqref{Non-DegenerateSolution}, we get the minimal-norm solution  
\begin{equation*}
	u_{0}^{\dagger}(t)=\frac{\pi}{2}H^{-1}_{11}h_{1}(t)=\sin t.
\end{equation*}

\section{Conclusion}
In this paper, we study the minimal-norm solution of the Fredholm integral equations of the first kind based on the \textit{H}-$H_{K}$ formulation. The basic idea of solving problems is to proceed from simple to difficult, specifically, from DKFIE to NDKFIE. For DKFIE, we obtain a closed-form representation of the minimal-norm solution using an operator method without conventional regularization methods or some other complicated methods. In addition, based on \textit{H}-$H_{K}$ formulation, we obtain the spatial structure of $N(L)^{\perp}$ and $R(L)$. Fortunately, the results that have been obtained can be extended easily to NDKFIE. Finally, multiple examples show that our proposed methods are feasible and effective.
\section*{Acknowledgements}
The research is supported by 2023 Guizhou University of Commerce  Project [grant number 2023ZKYB003],~Zhejiang Provincial Natural Science Foundation of China [grant number LQ23A010014].

\end{document}